\documentclass[11pt]{amsart}

\usepackage{verbatim}

\usepackage[backref,bookmarks=false]{hyperref}

\usepackage{color}

\numberwithin{equation}{section}

\newcommand{\ben}{\begin{enumerate}}
\newcommand{\een}{\end{enumerate}}

\newcommand{\bea}{\begin{eqnarray}}
\newcommand{\ba}{\begin{array}}
\newcommand{\bean}{\begin{eqnarray*}}
\newcommand{\ea}{\end{array}}
\newcommand{\eea}{\end{eqnarray}}
\newcommand{\eean}{\end{eqnarray*}}
\newcommand{\beq}{\begin{equation}}
\newcommand{\eeq}{\end{equation}}
\newcommand{\bthm}{\begin{thm}}
\newcommand{\ethm}{\end{thm}}
\newcommand{\blem}{\begin{lem}}
\newcommand{\elem}{\end{lem}}
\newcommand{\bprop}{\begin{prop}}
\newcommand{\eprop}{\end{prop}}
\newcommand{\bcor}{\begin{cor}}
\newcommand{\ecor}{\end{cor}}
\newcommand{\bdfn}{\begin{dfn}}
\newcommand{\edfn}{\end{dfn}}
\newcommand{\brem}{\begin{rem}}
\newcommand{\erem}{\end{rem}}
\newcommand{\bpf}{\begin{proof}}
\newcommand{\epf}{\end{proof}}
\newcommand{\bfact}{\begin{fact}}
\newcommand{\efact}{\end{fact}}
\newcommand{\bobs}{\begin{obs}}
\newcommand{\eobs}{\end{obs}}
\newcommand{\bexam}{\begin{exam}}
\newcommand{\eexam}{\end{exam}}
\newcommand{\bclaim}{\begin{claim}}
\newcommand{\eclaim}{\end{claim}}

\newtheorem{thm}{Theorem}[section]
\newtheorem{prop}[thm]{Proposition}
\newtheorem{lem}[thm]{Lemma}

\newtheorem{cor}[thm]{Corollary}
\newtheorem{dfn}[thm]{Definition}
\newtheorem{rem}[thm]{Remark}
\newtheorem{fact}[thm]{Fact}
\newtheorem{claim}[thm]{Claim}
\newtheorem{obs}[thm]{Observation}
\newtheorem{exam}[thm]{Example}

\newtheorem*{condition'}{Condition 2'}

 \newtheoremstyle{claimstyle}%
   {}
   {}
   {\normalfont}
   {}
   {\itshape}
   {.}
   { }
   {\thmnote{#3}}

\theoremstyle{claimstyle}

\alph{enumii} \roman{enumiii}

             \def\cB{\mathcal B}       
\def\cH{\mathcal H}                    
                   
                       \newcommand{\J}{\mathcal{J}}
              \def\cS{\mathcal S}             
                    
 \def\cI{\mathcal I}

                \def\Z{{\mathbb Z}}      
\def\C{{\mathbb C}}

\newcommand{\cbar}{\hat{{\mathbb C}} }

                \def\b{\beta}             
                          
                \def\Ga{\Gamma}           \def\l{\lambda} 
                         \def\Om{\Omega}
               \def\sg{\sigma}
                          
\def\ka{\kappa}

\newcommand{\ep}{\varepsilon}
\newcommand{\ph}{\varphi}
\newcommand{\al}{\alpha}
\newcommand{\ga}{\gamma}

\def\1{1\!\!1}

\def\and{\text{ and }}

\def\({\bigl(}                \def\){\bigr)}

                        \def\^{\tilde}

\def\ov{\overline}            \def\un{\underline}

\def\arg{\text{arg}}

\def\D{{\mathbb D}}

\def\${$ \displaystyle }

\newcommand{\jul}{\mathcal J}

\def\hypdim{HypDim(f)}

\def\Tract{\Omega}

\begin{document}

\title[Hyperbolic dimension of meromorphic functions]{A lower bound of the hyperbolic dimension for meromorphic functions having a logarithmic H\"older tract}


\author{Volker Mayer}
\address{Universit\'e de Lille, UFR de
  Math\'ematiques, UMR 8524 du CNRS, 59655 Villeneuve d'Ascq Cedex,
  France} \email{volker.mayer@univ-lille.fr \newline
  \hspace*{0.42cm} \it Web: \rm math.univ-lille1.fr/$\sim$mayer}

\date{\today} \subjclass{111}

\begin{abstract}  We improve existing lower bounds of the hyperbolic dimension for meromophic functions that have a logarithmic  tract $\Tract$ which is a H\"older domain. These bounds are given in terms of the fractal behavior, measured with integral means, of the boundary of $\Tract$ at infinity. 
 \end{abstract}

\maketitle

\section{Introduction}
Let $f:\C\to \cbar$ be a transcendental meromorphic function with Julia set  $\J_f$. The different fractal dimensions of such Julia sets has been intensively studied and the interested reader can consult the survey by Stallard \cite{Stallard-Survey}.
It is known that the Hausdorff dimension $Hdim(\jul_f)\geq 1$  for every entire function $f$ since then $\J_f$ contains non-degenerated continua \cite{Baker01}. Stallard showed that  $Hdim(\jul_f)>1$ for all functions of class $\cB$ which is the class of functions having a bounded set of singularities of the inverse. This result was later generalized by Rippon and Stallard to meromorphic functions having finitely many poles \cite{RS-2006} and by Bergweiler, Rippon and Stallard to functions having a logarithmic tract over infinity \cite{BRS-2008}.

The \emph{hyperbolic dimension}
of $f$, denoted $\hypdim$, is the supremum of the Hausdorff dimensions of all expanding repellers $X\subset \J_f$.
Given this definition, $\hypdim \leq Hdim (\jul_f )$ but still we have the same kind of result since 
Bara\'nski, Karpi\'nska and Zdunik \cite{BKZ-2009} have shown that 
$\hypdim >1$
for every meromorphic function $f$ having a logarithmic tract over infinity. 
This result is of remarkable generality.
Also it is sharp in the sense that $f_\l (z)=\l e^z$, $\l>0$ sufficiently small, has hyperbolic dimension
arbitrarily close to $1$, a result due to Karpinska \cite{karpinska-1999}. 

 Nevertheless, we improve this result still for functions having  a logarithmic tract $\Tract$ but with the additional assumption that one of these tracts has nice geometry, i.e. is a quasidisk, a John or, more generally, a H\"older domain (see Definition \ref{26}).

\bthm\label{thm main intro}
If the mermorphic map $f$ has a logarithmic tract over infinity, for example if $f\in \cB$, 
and if this tract is H\"older
then
$$\hypdim \geq \Theta \geq 1\,.$$
\ethm

But what precisely is this number $\Theta\,$? Let us for the moment simply mention that $\Theta\in [1,2[$ and that all values of this interval can occur. The precise definition is in Section \ref{15} and the final statement in Section \ref{sec 4}. 

\medskip

In order to obtain a good lower estimation of the hyperbolic dimension we take into account the geometry of a logarithmic tract. It appears that the size of the Julia set and of the hyperbolic dimension are influenced by the fractal behavior of the boundary of the tract near infinity. 
A logarithmic tract (over infinity) is a component $\Omega$ of $f^{-1}(\{|z|>R\})$ for some large $R$. The boundary $\partial \Om$
is a smooth curve. Nevertheless, seen from infinity $\partial \Om$ can have fractal behavior.
An appropriated way of measuring this  is to employ integral means and $\b$-numbers of natural rescalings of the Riemannian representation of the tract. This leads to a $\b$--number at infinity, called $\b_\infty$. Like 
in the classical case (see \cite{PommerenkeBook}  and \cite{Makarov98}), we will see that for H\"older tracts the function
$$t\mapsto \b_\infty (t)-t+1$$
has a unique zero $\Theta$ and this turns out to be a better lower bound for the hyperbolic dimension.

\medskip

Theorem \ref{thm main intro} is formulated for global functions $f$ since this is our main motivation. 
What really matters is the dynamics of $f$ in a logarithmic tract over infinity $\Omega$. So, we could
have considered in Theorem \ref{thm main intro} as well local functions like $F= f_{|\Om}$.  
Keeping the point of view of the local functions, we provide examples for which the bound $\Theta$ is expressed in terms of the Minkowski dimension of a quasicircle and such that the range of
this bound  is any number in $[1,2[$.

\medskip

We also investigate strict inequality between $HypDim(f)$ and $\Theta$ since this has important applications
in thermodynamical formalism (see \cite{MUpreprint3}). In general it is not clear whether this holds.
Bara\'nski, Karpi\'nska and Zdunik's do have strict inequality in their original result in \cite{BKZ-2009}.
This implies that strict inequality holds if $\Theta =1$. We generalize this in Proposition \ref{prop strict} to examples with general $\Theta$. The specificity of these examples is that they have tracts whose boundary have good additional geometric properties.


\section{Preliminaries}

This note concerns  transcendental meromorphic functions $f:\C\to \cbar$ having a logarithmic tract over infinity.
We denote $S(f)$ the singular set which is the closure of the critical and asymptotical values of $f$.
The precise definition of the Julia set, which will be denoted by  $\jul_f$, and other dynamical features of transcendental functions are contained in
 the survey by Bergweiler \cite{Bergweiler-survey}.
We use standard notations such as $\D(z,r)=\{|z]<r\}$, $A(r,R)= \D(z,R)\setminus \overline \D(z,r)$ and
$f\asymp g$ if $f/g$ are bounded below and above either by absolute constants or by constants where we indicate the dependance.

A key tool in dynamics is Koebe's Distortion Theorem. Theorem 1.3 in \cite{PommerenkeBook} is a version of it and from it follows that 
\beq\label{koebe}
\frac 18 (1-|z|) \leq \frac{|\psi '(z)|}{|\psi '(0)|}\leq \frac{2}{(1-|z|)^3} \quad \text{for every } \; z\in \D
\eeq
for every univalent map $\psi$ defined on the unit disk $\D$ with values in $\C$. 

An unbounded simply connected domain $\Tract \subset \C$ is a logarithmic tract over infinity of the function $f$ if there is 
$R>0$ such that $f_{|\Tract}:\Tract \to \{ |z| >R\}$ is a universal covering. 
Then there exists a conformal proper map $\tau :\Tract\to \{ \Re z > \log R\} $ such that $f=e^\tau$ on $\Tract$. 
Throughout the paper we will denote $\ph=\tau ^{-1}$ and suppose that $R=1$ so that 
$\ph : \cH=\{\Re z >0\} \to \Tract$. If $x>1$ and $z=x+iy\in \cH$, then it follows from \eqref{koebe} 
by considering the univalent map $\psi$ defined by $\psi (\xi)= \ph (x\xi + x+iy)$, $\xi\in \D$,
that
\beq\label{29}
\frac{1}{8x}\leq \frac{|\ph'(1+iy)|}{|\ph'(x+iy)|}\leq 2x^3\,.
\eeq

Every function of class $\cB$ has a logarithmic tract over infinity. 
Bishop's papers \cite{Bishop-EL-2015, Bishop-S-2016} give a complete new
understanding of this class of functions. He shows that every local model of a tract $\tau :\Tract\to \cH $
can be approximated by a global function $f\in \cB$ and even in $\cS$, the class of functions having only finitely many
singularities of the inverse. His papers complete former work by Rempe-Gillen \cite{Rempe-HypDim2}.
Contrary to the next definition, Bishop considers much more general models where $\Tract$ can be an arbitrary union of simply connected unbounded domains.

\bdfn\label{tract model}
A tract model $(\tau , \Tract )$ is a simply connected unbounded domain $\Tract$ along with
a conformal map $\tau : \Tract \to \cH$ such that $\tau^{-1}$ is proper. 
\edfn

This note considers H\"older tracts. Let us now explain what precisely this means. First of all,
a simply connected  hyperbolic domain $U\subset \C$ is H\"older if a conformal Riemann representation $\ph:\D\to U$
is H\"older:  there exists $H>0$ and $0<\al\leq 1$ such that 
\beq\label{holder 1}
|\ph (z_1)-\ph(z_2)|\leq H |\ph' (0)| |z_1-z_2|^\al \quad\text{for all}\quad z_1,z_2\in \D\,.
\eeq
It is a well known result by Hardy and Littlewood (see for example Theorem 5.1 in \cite{DurenHp}) that \eqref{holder 1} holds 
if and only if there exists $M>0$ such that
\beq\label{holder 2}
|\ph'(z)|\leq M|\ph' (0)| \left( \frac{1}{1-|z|}\right)^{1-\al} \quad \text{for all} \quad z\in \D
\eeq
and, moreover, the constants $M$ and $H$ depend mutually on each other. To be more precise,  \eqref{holder 1} implies  \eqref{holder 2} with $M=H$ and  \eqref{holder 2} implies  \eqref{holder 1} with $H= c\frac{M}{\al}$ where $c$ is an absolute constant.

In both conditions  \eqref{holder 1} and \eqref{holder 1} we have introduced the factor $|\ph' (0)|$ in order to make them scale invariant. If the image domain $U$ is replaced by $\l U$, $\l\neq 0$, then these two conditions still hold with unchanged constants $H$ and $M$. 

\blem\label{25}
If $\ph:\D\to U$ satisfies \eqref{holder 1} then
$$
 |\ph' (0)| \asymp dist  (\ph (0), \partial U)\asymp diam(U)
$$
where the involved constants depend only on $H$ and $ \al $.
\elem

\bpf
The relation $ |\ph' (0)| \asymp dist  (\ph (0), \partial U)$ results from Koebe's Distortion Theorem (see Corollary 1.4 in \cite{PommerenkeBook}).
One always has $dist(\ph (0), \partial U) \leq diam (U)$ and \eqref{holder 1} implies 
$$diam (U) \leq H |\ph'(0)| 2^\al\,.$$ Lemma \ref{25} follows.
\epf

A H\"older domain is simply the image of the unit disk by a H\"older map. But such domains are clearly bounded
whereas logarithmic tracts are unbounded domains.
One might try to use spherical geometry but this is not really adapted here.  A natural concept of 
H\"older tracts is the object of the following definition. It uses the notation 
$$Q_R=\{0<\Re z<2R \; , \; -R<\Im z < R\}\quad ,\quad R >0,$$ and a variant of \eqref{holder 1}: a conformal map $\ph:Q_1\to U$ is called $(H,\al )$--H\"older if
\beq\label{holder 1'}
|\ph (z_1)-\ph(z_2)|\leq H |\ph' (1)| |z_1-z_2|^\al \quad\text{for all}\quad z_1,z_2\in Q_1\,.
\eeq

\medskip

\bdfn\label{26}  
 A tract model $(\Tract , \tau )$ is  H\"older, more precisely $(H,\al )$--H\"older, if 
 for every $\kappa \in (0,1)$ there exists $m>0$ and $R_0\geq 1$ such that, for every $R\geq R_0$,
 \ben
 \item $\ph_R = \tau^{-1}\circ R :Q_1\to \Tract_R=\tau^{-1}( Q_R)$ satisifies  \eqref{holder 1'} and 
 \item $|\ph_R (z)|\geq m \, diam (\Tract_R )$ for every $z\in Q_1\setminus Q_{\kappa}$.
 \een
\edfn

\medskip

Item (1) simply means that $\Tract$ is exhausted by a family of uniformly H\"older domains $\Tract_R$
and item (2) is a quantitative version of the fact that $\ph$ is a proper map. Notice that it suffices to have 
item (2) for some $\kappa$, for example $\kappa=\frac 12$. Then one can easily recover the general case using \eqref{29} and Lemma \ref{25}.




\medskip

If $\ga$ is a rectifiable arc, then we denote $\ga (a,b)$ a subarc of $\ga$ joining $a,b\in \ga$ and by $|\ga (a,b)|$ its arc length. 
 A simply connected domain $U\subset \C$ is called $c$--John with base point $x_0\in U$, $c>0$, if for every $x_1\in U$ there exists a rectifiable arc $\ga\in U$ joining $x_0$ and $x_1$ such that 
$$dist(x,\partial U)\geq c \, |\ga (x , x_1) | \quad \text{for every} \quad x\in U\,.$$
John domains are H\"older but not vice versa (see \cite{BP82} and \cite{PommerenkeBook}). We can therefore also consider the following class of John tracts:

\bdfn\label{26 John}  
 A tract model $(\Tract , \tau )$ is  $c$--John if for every $R\geq 1$ the domain 
 $ \Tract_R=\tau^{-1}( Q_R)$ is $c$--John with base point $z_R= \tau^{-1}( R)$.
 \edfn

\medskip

A tract $\Tract $ of $f$ is a quasidisk and $\Gamma =\partial \Tract$ a quasicircle if  the conformal map $\tau :\Tract\to \cH$ has a quasiconformal extension to the plane. Quasidisks are John domains but the boundary of John domains is not necessarily a Jordan curve. 
Logarithmic tracts that are quasidisks are John tracts and thus H\"older tracts. So, logarithmic tracts that are quasidisks are examples of tracts that fit into our context. Here are more examples that the reader may have in mind.

\bexam\label{CE}
Let $P_c(z)=z^2+c$ be a polynomial with connected Julia set and having a repelling fixed point $p\in \partial A(\infty )$,
 the attracting bassin of infinity of $P_c$. Then it is known that
\ben
\item $A(\infty )$ is a John domain if and only if $P_c$ is semi-hyperbolic, a result by Carleson, Jones and Yoccoz \cite{CJY94}, and
\item if $p_c$ is Collet-Eckmann, then $A(\infty )$ is a H\"older domain as Graczyk and Smirnov showed in  \cite{GS98}
\een
There exists a linearizing entire function 
$f$ such that $f\circ \l = P_c\circ f$ where $\l$ is the multiplier of $P_c$ at $p$. If (1) or (2) holds, then 
$\Tract =f^{-1} (A(\infty ))$ is respectively a John or H\"older tract (\cite{MUpreprint3}
contains all the details and further informations on that).

 Here $\partial \Tract$ is a fractal and $\Tract $ is not exactly a tract as defined above. So, let $R>0$ sufficiently large such that
 $\{|z|=R\} \subset A(\infty )$ Then, the Julia set of $P_c$ being connected, $\{|z|>R\} \subset A(\infty )$
 and $\Tract = f^{-1}(\{|z|>R\} )$ is a tract as defined above. Here $\partial \Tract$ is smooth. But let us consider 
 $$
 \partial \Tract \cap \left\{ |\l^n|\leq |z| < |\l^{n+1}|\right\}\quad , \quad n\geq 1\,.
 $$
If we rescale this set to unit size by multiplying it with $\l^{-n}$ then it follows from the linearizing functional equation 
that this set ressembles more and more the fractal boundary of $A(\infty )$. In fact, these rescalings behave exactly 
like equipotential lines.
\eexam

\section{Integral means at infinity}\label{15}

Let $\ph:\cH =\{\Re z >0\} \to \Tract$ be a conformal map with $\ph(\infty ) =\infty$, i.e. such that $\lim_{z\to \infty} |\ph(z)|=\infty$. In practice and up to normalization by an additive constant $\ph$ will be the reciprocal of the conformal map $\tau$ used in the above definition of the logarithmic tract. Let $\kappa ^{-1}= 1+2\pi $, let   $r>0 $, $t\geq 0$  and consider
\beq \label{1}
I_\ph (r,t) = \int _{\ka}^1 |\ph '(r+iy)|^t dy \;\; \text{ and }\;\; \b_\ph (r,t) = \frac{\log I_\ph (r,t) }{\log 1/r}\,.
\eeq

\brem\label{17} The classical setting (see  \cite{Makarov98} and \cite{PommerenkeBook})
considers conformal maps on the unit disk and integral means on the circles of radius $0<r<1$.
Then the behavior of the classical $\b$--function
\beq\label{beta}
\b (t)=\b_\ph (t)=\limsup_{r\to 1}\b_\ph (r,t)
\eeq
 is still not completely clear. If the range is a H\"older domain, then one knows
 (see \cite{Makarov98} and page 241 as well as Corollary 10.18 in \cite{PommerenkeBook})  that
  \ben
\item[a)]  there exists $0\leq k <1$ such that $\b (t+s)\leq \b (t) +k s$ which implies that
\item[b)] $t\mapsto \b (t)-t+1$ is strictly decreasing and, consequently, that there exists a unique $\Theta$ such that $\b (\Theta )= \Theta-1$.
\een
Moreover,  the number $\Theta$ is the Minkowski dimension of the boundary of the image domain.
\erem

Since we are interested in the fractal behavior of $\Tract$ near infinity it is natural to consider rescalings.
Let $R\geq R_0$ and consider
\beq\label{4}
g_R := \frac{1}{diam(\Tract _R )}\circ   \ph \circ R : Q_1 \to \tilde \Tract_R= \frac{\Tract_R}{diam(\Tract _R )}\,.
\eeq
Since the tract is H\"older, this map satisfies  
\beq\label{4'}
|g_R(z)|\geq m \quad \text{for every } \quad z\in Q_1\setminus Q_\kappa
\eeq
and \eqref{holder 1'} with $|g_R'(1)|\asymp 1$ and with constant independent of $R\geq R_0$.
Let $z_1,z_2\in Q_1\setminus Q_\kappa$. Then
$|g_R(z_1)-g_R(z_2)|\leq \tilde H$ with $\tilde H = \tilde H (H)$.  Therefore,
\beq\label{22}
\frac{|\ph (Rz_1)|}{|\ph (Rz_2)|} \leq  1+ \frac{\tilde H}{|g_R(z_2)|}\leq 1+ \frac Hm =: M
\quad \text{for every} \quad z_1,z_2\in Q_1\setminus Q_\kappa\,.
\eeq

\medskip

If the map $g_R$ is considered as a map on the segment $\{\Re z =1/R\}\cap Q_1$ it is a rescaling
of the map $\ph$ defined on $\{\Re z = 1\}\cap Q_R$. We now apply the definitions in \eqref{1} to these rescalings
and set  $I^+_R(t):= I_{g_R}(1/R,t)$. The  symmetric analog of this is
$$
I^-_R(t):=\int _{-1}^{-\ka} |g_R '(1/R+iy)|^t dy\,.
$$
Consider
\beq \label{2}
I_R(t)=\max \big\{I^-_R(t), I^+_R(t)\big\}
 \quad \text{and }\; \b_{R}(t)=\frac{\log I_R(t)}{\log R}\,.
 \eeq
and define
\beq \label{3}
\b_{\infty}(t)= \limsup_{R\to \infty}  \b_{R}(t)\quad , \; \;t\geq 0\,.
\eeq

In these definitions we did two arbitrary choices. First of all, 
we fixed  $\kappa^{-1} = 1+2\pi$  but any $\kappa \in (0,1)$ can be used here 
and will lead to the same number $\b_\infty$.
The second choice is that
we considered the map $\ph$ on the line $\{ \Re z =1\}$. 
Although we will use it differently,
the following Lemma explains that one could 
replace this line by any other line $\{ \Re z =s\}$ without changing $\b_{\infty}$.

\blem
\label{6}
For every $x>1$, $t\geq0$, $\kappa \in (0,1)$ and $T\geq 1$ we have that
$$
2^{-t} x^{-3t}\leq \frac{\int _{\kappa T}^T |\ph'(x+iy)|^tdy }{ \int _{\kappa T}^T |\ph'(1+iy)|^tdy}\leq 8^tx^{t}\,.
$$
\elem

\bpf
Follows directly from \eqref{29}.
\epf

The function $\b_\infty$ has exactly the same good properties than the standard $\b$-function if $\Tract$ is a H\"older tract.
\bprop\label{prop 3}
Let $\Tract$ be a H\"older tract. Then items a) and b) of Remark \ref{17} hold with $\b$ replaced by $\b_\infty$
and there exists a unique number $\Theta$ such that $$\b_\infty (\Theta ) -\Theta +1=0\,.$$
\eprop

\bpf
 This mainly follows from  \eqref{holder 2} by following the argument in  \cite{PommerenkeBook}.
 Here is a way to adapt this proof to the present setting.
 
 Let $R >1$ and let $g_R$ be the rescaled map defined in \eqref{4}.
 We can consider $g$ as a map on the larger square $Q_{1+1/R}$ and then it still is H\"older and we may assume that
 the constants $H,\al$ are unchanged. The segment
 $\big\{ \Re z =1/R\big\}\cap Q_1$ has distance $1/R$ to the boundary of $Q_{1+1/R}$.
 On the other hand, the Riemann map from $\D$ onto $Q_{1+1/R}$ mapping $0$ to $1$ is, uniformly in $R>1$, bi-H\"older.
 It therefore follows from  \eqref{holder 2} that there exists $\tilde M>0$ such that
 $$
 |g_R'(1/R+iy)|\leq \tilde M R^{1-\al }\quad \text{for every $|y|\leq 1$ and every $R>1$.}
 $$
 This inequality allows the following estimation: if $t,s>0$, then
$$
I_R^+(t+s)= \int _{\kappa}^1 |g_R '(1/R+iy)|^{t+s} dy\leq \tilde M^s R^{(1-\al)s} I^+_{R}(t)
$$
and the same inequality is true for $I_R^-$ hence also for $I_R$.
Consequently,
$$\b_\infty (t+s)\leq \b_\infty (t) +({1-\al }) s= \b_\infty (t) +k s$$ where $k=1-\al\in [0,1[$. The proof of the remaining affirmations of the Proposition \ref{prop 3} is now straightforward.
\epf


\section{The result and examples}\label{sec 4}

The aim of this note is to establish the following result. Inhere we associate to a logarithmic tract 
the function $\b_\infty$ as explained in the previous section.

\bthm\label{thm main}
Let $f:\C\to\cbar$ be a meromorphic function having a logarithmic H\"older tract $\Tract$ over infinity. Let $\Theta$ be the unique zero of the function $t\mapsto \b_\infty (t)-t+1$. Then,
$$\hypdim \geq \Theta \geq 1\,.$$
\ethm

\smallskip
As already explained in the Introduction, this result also holds for local functions $F=e^\tau$ where $(\Om , \tau )$
is a tract model. The proof of Theorem \ref{thm main} does indeed only us the behavior of $f$ in the tract.

In many situations the $\Theta$--number can be determined. For H\"older domains $\Theta$ equals the Minkowski dimension of the boundary of the domain (see Remark \ref{17}).
The following example is similar to the linearizing functions in Example \ref{CE}.
 
 \bexam \label{example}
 Let $\sg$ be a bounded quasicircle such that $\sg \cap [0,\infty[  =\{1\}$ 
 and such that the origin is in the bounded component $U$ of $\C\setminus \sg$.
 Let $\b_\sg$ the $\b$--function of the interior of this curves as defined in \eqref{beta}. More precisely, let $\ph:\D\to U$ be any conformal map and take then the associated $\b$--function.
Consider then $\ga_1$ the lift by the exponential map of $\sg$
 having endpoints $2i\pi$ and $4i\pi$. Define $\ga_k=2^{k-1} \ga_1$, $k\in \Z$, and
 $$\Gamma ^+ =\bigcup_{k\in \Z} \ga _k \quad \text{and} \quad \Ga = -\Ga^+\cup \{0\}\cup \Ga^+\,.$$
 This curve $\Ga$ is a quasicircle. If $\Tract$ is one of the connected components of $\C\setminus \Ga$
 and $\tau: \Tract \to \cH$ a conformal map then $(\tau, \Om )$ is a tract model. By construction of $\Ga$, a simple verification shows that the $\b_\infty$--function of $(\tau, \Om )$ is
 equal to $\b_\sg$. Theorem \ref{thm main} and Corollary 10.18 in \cite{PommerenkeBook} imply now
 $$HypDim(f) \geq \Theta_f = Mdim ( \sg )\,.$$
\eexam

These examples  can be easily generalized in many ways simply by taking for every $k\in \Z$
a different generating curve $\sg_k$ in order to define $\ga_k$. 

\medskip

 Contrary to the result in \cite{BKZ-2009} by Bara\'nski, Karpi\'nska and Zdunik, it is not clear whether strict inequality holds in Theorem \ref{thm main}. The case $\Theta =1$ is special since the boundary of the tract has infinite length. 
Let us add a similar assumption for the generating curve $\sg$ in the preceding example.

\bprop \label{prop strict}
Suppose $f$ is like in the Example \ref{example}, that $$\Theta = Mdim(\sg ) <2$$ and suppose that
$\sg$ has $\Theta$--Hausdorff measure $HM^\Theta (\sg )>0$. Then
$$\hypdim > \Theta\,.$$
\eprop

\medskip

 \section{Proof of Theorem \ref{thm main}}

The main step in this proof is to generate a good iterated function system.
Most often this is done after a logarithmic change of coordinates. We avoid this here but use it later implicitly by employing the  metric $|dz|/|z|$. 

The setting is the following. By assumption, $f$ has a logarithmic H\"older tract $\Tract$. Let again $\ph$ be the inverse of 
$\tau=\log f$ (modulo normalization by an additive constant) restricted to $\Tract$. For $s>0$, we denote
\beq \label{9}
\Gamma_s =\ph (\{ \Re z = s\}) \; , \; \Gamma_s^- =\ph (\{ \Re z = s\; , \;  \Im z <0\} )
\; \text{ and } \; \Gamma_s^+ =\Gamma_s \setminus \Gamma_s^-\;.
\eeq

The first step is of technical nature and it is similar to Corollary 3 in \cite{BKZ-2009}. 

\blem\label{8}
There exists $R_0\geq 1$ such that for every $R\geq R_0$
\beq\label{33.33}
\Gamma_{\log R}^-\cap \{|z|=R\}\neq\emptyset\; \text{ and }\; \Gamma_{\log R}^+\cap \{|z|=R\}\neq\emptyset\,.
\eeq
Moreover, increasing $R_0$ if necessary, there exists $C<\infty$ independent of $R\geq R_0$  such that $R\leq C (\Im a_0)^2$
for every $a_0 = \log R + i \Im a_0$ with $$b_0=\ph (a_0)\in \Gamma_{\log R}^+\cap \{|z|=R\}\,.$$
\elem

\bpf
Let $T\geq 1$ and remember that $\ph_T=\ph \circ T : Q_1\to \Tract_T$ satisfies the H\"older inequality \eqref{holder 1'}:
$
|\ph_T(\xi_1) -\ph_T(\xi_2)|\leq H |\ph'_T (1)| |\xi_1-\xi_2|^\al$, $\xi_1,\xi_2\in Q_1$.
Clearly $|\ph_T '(1)|= T|\ph '(T)|$ and \eqref{29} implies that 
$
|\ph '(T)|\leq 8 |\ph '(1)|T
$.
Therefore,
\beq\label{33.3}
|\ph(z_1)-\ph(z_2)|\leq  8H |\ph '(1)| T^{2-\al} |z_1-z_2|^\al \; \; , \quad z_1,z_2\in Q_T\,.
\eeq

Now, let $R>e$ and set $T=\log R$. Inequality \eqref{33.3} implies
$$
|\ph (T)-\ph (1)|\leq 8H |\ph '(1)| T^{2-\al} |T-1|^\al\leq 16H |\ph '(1)| T^{2} 
$$
which shows that $|\ph(\log R )| \leq |\ph (1)| + 16H |\ph '(1)| (\log R)^2$. 
It suffices now to choose $R_0>1$ such that 
$|\ph (1)| + 16H |\ph '(1)| (\log R)^2 < R$ for all $R\geq R_0$
since then $|\ph (\log R)| <R$ which  implies \eqref{33.33} the map $\tau =\ph^{-1}$ being a proper map.

Consider now a point $b_0=\ph (a_0)\in \Gamma_{\log R}^+\cap \{|z|=R\}$ and set
$$
T=2 \max \{\log R \, , \, y\} \text{ where } \log R = \Re a_0 \text{ and } y=\Im a_0 \geq 0 \,.
$$
Then $a_0\in Q_T$ and we can use again \eqref{33.3} with $z_1=a_0$ and $z_2=0$:
$$
|\ph (a_0)-\ph (0)|\leq 8H |\ph '(1)| T^{2-\al}|a_0|^\al\leq 16 H |\ph '(1)| T^2\,.
$$
Since $R=|b_0|=|\ph (a_0)|$ it follows that 
$$R\preceq \max \{ \log R \; ,\,  y\} ^2\,.$$
Thus $R\preceq y^2= (\Im a_0)^2$ provided $R$ is sufficiently large.
\epf

In the following we assume that $R\geq R_0$ with $R_0$ such that both Lemma \ref{8} 
and the Definition \ref{26} apply with $\kappa^{-1} = 1+ 2\pi$.
Consider a point
$$b_0\in \Gamma_{\log R}^+\cap \{|z|=R\} \text{ and set } a_0=\ph^{-1}(b_0)\,.$$
Here we used the positive part of $\Gamma_{\log R}$. This can be done if the maximum
 in \eqref{2} is atteint by $I_R^+$. We may and will assume that this is the case
 since otherwise it suffices to use here $\Gamma_{\log R}^-$ instead of $\Gamma_{\log R}^+$.

Define this times $T=\kappa^{-1} \Im a_0$. Let also $N$ be the maximal integer such that $ \Im a_0+2\pi (N+1) \leq T$
 and define the points
 \beq\label{a_ks}
 a_k=a_0+ 2i\pi k \; \text{ and }\;b_k=\ph (a_k) \; \text{ , }\;
 k\in \{ 0, 1,..., N\}\,.
 \eeq
All the points $b_k$ have the same image 
$$\tilde b = f(b_k)\; , \;\;k\in \{ 0, 1,..., N\}\,.$$
Since $a_0,...,a_N\in Q_T\setminus Q_{\kappa T}$ and since $|b_0|=R$, \eqref{22} implies that there exists $M\geq 1$
independent of $R$ such that all the points
$b_k$ are in the annulus $A(R/M , MR )=\{\frac{1}{MR}<|z|<MR\}$. Let
$$ A = A(R/2M , 2MR ) \,.$$
Let $l\in \{0, ..., 7\}$, the precise value of this integer will be chosen later.  Associated to it is the set $E$ of $k\in \{0, ..., N\}$ such that $arg (b_k) \notin [l \frac \pi 4 , (l+1) \frac \pi 4[$.
Consider the ray
 $\gamma = \gamma_l = \{z\in A \; , \arg z = (2l+1)\pi /8\}$ and the domain
$$V = A \setminus \gamma \; .$$
We may assume that $\tilde b \notin \gamma$ since otherwise it suffices to replace $\ga$ by an other ray 
that is arbitrarily close to it. Finally, let $U_k$ be the connected component of $e^{-1} (V)$ that contains $a_k$ and 
$V_k=\ph (U_k)$, $k\in E$. Clearly, the whole construction does depend on $R>R_0$.

\blem\label{30}
There exists $R_1>R_0$ such that for every $R>R_1$ we have $$\overline V_k\subset V \quad , \quad  k\in E\,.$$
\elem

\bpf
Let $k\in E$. First of all, $diam(U_k)\asymp \log (2M)$ and clearly $U_k\subset Q_T$ if $R$ is sufficiently large
since then $T$ is large by Lemma \ref{8}. It precisely shows that  
\beq\label{34}
R\preceq T^2\,.
\eeq
This inequality also shows that $\Re a_0 =\log R <  \Im a_0 = \ka T$ and thus 
$a_0 \in Q_T\setminus Q_{\kappa T}$.
Since $\kappa <\frac12$, $T\in  Q_T\setminus Q_{\kappa T}$. Therefore,  \eqref{22}
applies and yields $\frac 1M \leq \frac{|\ph (T)|}{|\ph (a_0)|}\leq M$ which itself implies $\frac RM \leq |\ph (T)| \leq MR$. Consequently,
$$
dist( \ph (T) , \partial \Tract_T )\leq |\ph (T)-\ph (0)| \asymp R\,.
$$
By Lemma \ref{25}, $dist( \ph (T) , \partial \Tract_T )=dist( \ph_T (1) , \partial \Tract_T )\asymp |\ph_T'(1)|$ and thus
$$
 |\ph_T'(1)| \preceq |\ph (T)-\ph (0)| \asymp R\,.
$$
On the other hand, the H\"older property (1) of Definition \ref{26}
 implies that
$$
|\ph (z_1)-\ph(z_2)|\leq H |\ph_T'(1)|\left( \frac{|z_1-z_2|}{T}\right)^\al
$$
and thus
$$
|\ph (z_1)-\ph(z_2)|\preceq \frac{R}{T^\al} |z_1-z_2|^\al \quad , \quad z_1,z_2\in Q_T\,.
$$
This relation holds in particular for all $z_1,z_2\in U_k$ which implies
$$
diam (V_k)= diam(\ph (U_k)) \preceq  \frac{R}{T^\al}\,.
$$
Therefore, $diam (V_k)/R\to 0$ as $R\to \infty$.
This completes the proof since $b_k\in V_k$ and 
$dist (b_k, \partial V)\geq c R$ with $c$ some absolute constant. 
\epf

Let in the following $R>R_1$. The map $f:V_k\to V$ is conformal. Denote $\Phi_k:V\to V_k$ the inverse map. Since the ratio between the inner and outer radius of 
$A$ does not depend on $R$, the maps $\Phi_k$ have bounded distortion and this distortion is independent of $R> R_1$.
All in all we have:

\bprop\label{13}
The maps  $\Phi_k:V\to V_k$, $k\in E$, define a conformal iterated function system as defined in \cite{MauldinUrb02} with limit set a subset
of the Julia set of $f$. 
\eprop

It remains to estimate the Hausdorff dimension of the limit set of this iterated function system.
In order to do so we can work with the cylindrical metric $|dz|/|z|$. The derivative $|f'(z)|_1$ of $f$ 
at a point $z\in \Tract$ with respect to this metric is
$$
|f'(z)|_1=\frac{|f'(z)|}{|f(z)|}|z|= \frac{|e^{\tau (z)}\tau'(z)|}{|e^{\tau(z)|}}|z|
=\frac{|\ph (w)|}{|\ph'(w)|}\;  \text{ with } \; w=\tau(z)
\,.
$$
Let $t>0$. We have to estimate $\Sigma_{t,l}=\sum_{k\in E} \big|\Phi_k'(w)\big|_1^t $ , $w\in V$.
Since the system has uniformly in $R>R_1$ bounded distortion it suffices to consider an arbitrary point $w\in V$.
We will take $w=\tilde b$ and take the notations of \eqref{a_ks}. In particular, $\Phi_k(\tilde b)= b_k$  if $k\in E$
and we have to consider$$
\Sigma_{t,l}=\sum_{k\in E} \big|\Phi_k'(\tilde b)\big|_1^t =\sum_{k\in E}\left| \frac{\ph' (a_k)}{\ph (a_k)}\right|^t
\,.
$$
Remember that $l$ is any integer between $0$ and $7$. We now suppose that $l$ is chosen such that 
$$
 \frac 87\Sigma_{t,l} \geq \Sigma_t = \sum_{k=0}^N\left| \frac{\ph' (a_k)}{\ph (a_k)}\right|^t\,.
$$
By  Koebe's Distortion Theorem $|\ph'(a_k)|\asymp |\ph'(a_k+iy)|$ for $0\leq y\leq 2\pi$. 
Item (2) of Definition \ref{26} and $a_k\in Q_T\setminus Q_{\kappa T}$ imply $|\ph (a_k)| \asymp diam(\Tract_T)$.
Consequently,
$$
 \Sigma_t \asymp \frac {1}{diam(\Tract_T)^t} \int _{\kappa T}^T |\ph '(\log R +iy)|^t dy\,.
$$
Lemma \ref{6} yields
$$
\frac{1}{(2(\log R)^3)^t}\int _{\kappa T}^T |\ph '(1 +iy)|^t dy\leq
\int _{\kappa T}^T |\ph '(\log R +iy)|^t dy\,.
$$
So, let us consider
$$
 \frac {1}{diam(\Tract_T)^t} \int _{\kappa T}^T |\ph '(1 +iy)|^t dy= 
 T^{1-t}\int _\kappa^1 \left|g_T '\left( \frac 1T +iy\right)\right|^t dy
$$
where the rescaled function $g_T$ has been defined in \eqref{4}.
We assumed that the maximum in \eqref{2} is atteint by $I_T^+$ and thus
$$
 \frac {1}{diam(\Tract_T)^t} \int _{\kappa T}^T |\ph '(1 +iy)|^t dy= T^{1-t}I_T^+(t)= T^{1-t+\beta_T (t)}
$$
from which follows that
\beq\label{14}
\frac{1}{(\log R)^{3t}} T^{1-t+\beta_T (t)}
\preceq \Sigma_t \,.
\eeq
Let $t<\Theta$. Then there exists $\ep >0$ and $T>1$ arbitrarily large  such that $$\b_T(t)-t+1 \geq \ep>0\,.$$ 
 Combining \eqref{14} and \eqref{34} shows that for every $A>1$ there exists
 $R$ and thus $T$ sufficiently large such that
$$\Sigma_{t,l}\geq \frac 78 \Sigma_{t}\geq A \,.$$
This implies that the iterated function system for this choice of $t,l$ has positive pressure 
and thus the Hausdorff dimension of the limit set of this system is strictly larger than $t$
(see \cite{MauldinUrb02} or \cite{PUbook}). We therefore proved that
$HypDim(f)\geq \Theta$.

\

It remains to show that $\Theta \geq 1$. Let $T_n=\kappa^{-n}$ and, for simplicity, let 
$g_n$ be the rescaled map $\frac{1}{diam(\Tract _{T_n} )}\circ \ph \circ T_n$. 
Remember that $|g_n(z)|\geq m$ for all $z\in Q_1\setminus Q_\kappa$ (see \eqref{4'}) and that  $|g_n'(1)|\asymp 1$.
We thus may assume that the maps $g_n$ satisfy the H\"older inequality  \eqref{holder 1'} without the factor $|g_n'(1)|$.

Therefore, $|g_n(\kappa^n+i)|\geq m$ and $|g_n((\kappa^n+i)\kappa ^N ) | \leq 2H\kappa ^{N\alpha } + |g_n(0) | $. Since $\lim_{n\to\infty}|g_n(0) | =0$ there exists $N$ and $n_0\geq 2N$ such that $|g_n((\kappa^n+i)\kappa ^N ) | \leq \frac m2$ for every $n\geq n_0$. Consequently,
$$|g_n(\kappa^n+i)-g_n((\kappa^n+i)\kappa ^N ) | \geq \frac m2 \quad\text{ for every }\quad n\geq n_0\,.$$
This implies that we can choose for every $n\geq n_0$ an integer $l\in \{0,...,N\}$ such that 
$|g_n((\kappa^n+i)\kappa^l)-g_n((\kappa^n+i)\kappa ^{l+1} ) | \geq \frac m{2N}$ from which we get
$$
\int _{\kappa^{l+1}}^{\kappa^l} |g_n'((\kappa^n+i)y)|dy \geq \frac m{4N}:=m'>0\,.
$$
Since $l\leq N$ and since $n_0\geq 2N$, 
it follows from Lemma \ref{25} and Koebe's Distortion Theorem that $$diam(\Tract_{T_{n-l}})\asymp diam (\Tract_{T_n})$$ with constants independent of $n\geq n_0$. A simple change of variable and once more bounded distortion gives
$$
I^+_{\kappa_{l-n}}(1)=
\int _{\kappa}^1 |g_{n-l}'(\kappa^{n-l}+iy)|dy 
\asymp\int _{\kappa^{l+1}}^{\kappa^l} |g_n'((\kappa^n+i)y)|dy \geq m'\,.
$$
We showed that there is $m''>0$ and arbitrarily large integers $\nu$ such that 
$$I_{\kappa^{\nu}}(1)\geq I^+_{\kappa^{\nu}}(1) \geq m'' >0\,.$$
Consequently, $\b_{T_\nu}\geq \frac {\log(m'')}{\log T_\nu }$ for all these integers $\nu$ and thus $\b_\infty (1)\geq 0$.

 \section{Proof of Proposition \ref{prop strict}}
Since we deal with a particular function it is possible to construct an other iterated function system
which then allows to conclude. 

\medskip

We use the setting and the notations from Example \ref{example} . The domain $\Om$ is one of the components of the complement of $\Ga$.
By construction, this curve and one of the domains of $\C\setminus\Ga$ does not intersect the ray $]-\infty , 0[$. We thus may assume that 
$]-\infty , 0[\; \cap \;\Om =\emptyset$.
For $1\leq \un S<\ov S$, define the domain
$$U=U(\un S, \ov S )=\Big\{ \un S<|z|<\ov S\Big\}\setminus ]-\ov S , -\un S[\,.$$
We lift this domain with different inverse branches of $F=e^\tau$ and choose among them thoses that stay inside $U$.
First, $$V_0=\{\log \un S<\Re z <\log \ov S \; , \; -\pi <\Im z< \pi\}$$ and $V_k= V_0+2ki\pi$, $k\in \Z$, are the components of 
$e^{-1}(U)$. Let $\ph: \cH\to \Om$ be again the inverse of the conformal map $\tau$. We may assume that $\ph (0)=0$. The domains we look for are
$$  U_k = \ph (V_k) \;\; , \;\; k\in \Z\,.$$
In order to select among the domains $U_k$ sufficiently many that stay in $U$ we have to discuss $\ph$,
make a choice for $\un S, \ov S$ and select then domains $U_k\subset U$.

The domain $\Om$ has the  invariance
$2\Om = \frac12 \Om=\Om$ and the points $b_k=2^ki\pi$, $k\in \Z$, are in the boundary $\Ga$.
Let $a_k=\ph^{-1}(b_k)$, $k\in \Z$, and $\mu= a_2/a_1$. 
This invariance implies that the conformal map $\ph_1=\frac12 \ph \circ \mu$ also maps $\cH$ onto $\Om$ and thus $\ph^{-1}\circ \ph_1$
is a conformal map of $\cH$ that fixes the three points $0, a_1$ and $\infty$. Consequently $\ph_1=\ph$ and it follows inductively that
\beq\label{functional equation 1}
\frac{1}{2^n}\ph \mu^n =\ph \quad \text{and} \quad a_n = \mu^{n-1}a_1\;\text{ for every }\; n\geq 1\,.\eeq
We may assume that $a_1=i$ so that $a_n =\mu^{n-1}i$.

Fix $0<\un S_0<\ov S_0$ such that 
\beq\label{functional equation 2}
 \ph (Q_{2\mu} \setminus Q_{1/\mu}) \subset \D_{\ov S_0-1}\setminus 
\ov\D _{\un S_0}
\,.
\eeq
Take then $\un S , \, \ov S$ of the form $\un S = 2^N \un S_0$, $\ov S= 2^{2N}\ov S_0$, $N>1$, and let
$$k\in  \cI= \bigcup_{n=N}^{2N-1}\cI_n= \bigcup_{n=N}^{2N-1} \frac{1}{2\pi}\big[\mu^{n}, \mu^{n+1}\big[\,.$$
Let in the following $N$ be sufficiently large such that 
\beq\label{functiona 2} 
\mu^{-N}\log \ov S < s_0
\eeq 
where $s_0>0$ is some small number that will be determined  later on.
Then, 
 for every $k\in \cI_n$, $n\in \{N,..., 2N-1\}$, $\mu^{-n} V_k \subset  Q_{2\mu} \setminus Q_{1/\mu}$ and thus, because of
 \eqref{functional equation 1} and  \eqref{functional equation 2},
$$U_k=\ph (V_k) = 2^n \ph (\mu^{-n} (V_k))\subset 2^n \ph (Q_{2\mu} \setminus Q_{1/\mu})\subset U\,.$$
Let $\psi_k:U\to U_k$ be the inverse of $F_{|U_k}:U_k\to U$. Then, the collections of maps $\psi_k$, 
$k\in \cI$,
defines again an iterated function system. 
The next result completes the proof of Proposition \ref{prop strict}

\blem
The Hausdorff dimension of the limit set of $(\psi_k)_{k\in \cI}$ is strictly larger than $\Theta$.
\elem

\bpf
We have to estimate 
$$
\sum_{k\in \cI} |\psi_k'(w)|^\Theta_1 = \sum_{k\in \cI} |F'(z_k )|^{-\Theta}_1 = 
\sum_{k\in \cI} \left(\frac{|\ph'(\xi_k )|}{|\ph (\xi _k)|)}\right)^{\Theta}
$$
where $w\in U$, $\xi_k\in e^{-1}(w_0) \cap V_k$ and $z_k=\ph (\xi _k)$. 

Let $n\in \{ N,..., 2N-1\}$.
If $k\in \cI_n$ then $|\ph (\xi _k)| \asymp 2^n$. Thus
$$
\sum_{k\in \cI_n} \left(\frac{|\ph'(\xi_k )|}{|\ph (\xi _k)|)}\right)^{\Theta}
\asymp \frac{1}{2^{n\Theta}}\sum_{k\in \cI_n} |\ph'(\xi_k )|^\Theta\asymp 
\frac{1}{2^{n\Theta}} \int _{\mu^n}^{\mu^{n+1}} |\ph ' (\log|w| +iy)|^\Theta dy.
$$
Using \eqref{functional equation 1} and a change of variable it follows that
$$
\sum_{k\in \cI_n} \left(\frac{|\ph'(\xi_k )|}{|\ph (\xi _k)|)}\right)^{\Theta}
\asymp \mu^{(1-\Theta)n} \int _1^\mu |\ph'(s+iy)|^\Theta dy
$$
where $s= \log |w| \mu^{-n}$. On the other hand, 
$$
 \int _1^\mu |\ph'(s+iy)|^\Theta dy \asymp s \sum_{k=0}^{[\frac{\mu-1}{s}]} |\ph'(s+i(1+sk))|^\Theta\,.
$$
The curve $\Gamma$ being a quasicircle, the later sum can be associated with $\sum_k diam (D_k)^\Theta$
where the $D_k$ are disks that cover $\gamma_1$ and have $diam(D_k)\asymp  |\ph'(s+i(1+sk))| s$.
By assumption, $HM^\Theta (\gamma_1)=m>0$ which implies that there exists $r_0>0$ such that 
$\sum_k diam (D_k)^\Theta \geq m/2$ provided $diam(D_k)\leq r_0$, $k\in \{0, ..., [\frac{\mu-1}{s}]\}$.
Uniform continuity of $\ph$ in a compact neighborhood in $\ov\cH$ of $[i,i\mu]$ allows to choose $s_0>0$
such that the disks $D_k$ above have $diam(D_k)<r_0$ if $s<s_0$. Assume now that $N$ is chosen such that 
\eqref{functiona 2} holds with this number $s_0$. Then $ s= \log |w| \mu^{-n} <s_0$ for all $w\in U$ and 
$n\in \{N,..., 2N-1\}$. Therefore,
$$
\sum_{k\in \cI_n} \left(\frac{|\ph'(\xi_k )|}{|\ph (\xi _k|)}\right)^{\Theta}
\asymp  \mu^{(1-\Theta)n} s^{1-\Theta}\sum_k diam (D_k)^\Theta \geq \frac{m}{2}  (\log |w|)^{1-\Theta }\asymp N^{1-\Theta }
$$
from which follows that 
$$
\sum_{k\in \cI} |\psi_k'(w)|^\Theta_1 =\sum_{n=N}^{2N-1}\sum_{k\in \cI_n} \left(\frac{|\ph'(\xi_k )|}{|\ph (\xi _k|)}\right)^{\Theta}
\succeq  N^{2-\Theta}\,.
$$
Since $\Theta <2$ this shows that we can choose $N$ sufficiently large such that 
$$
\sum_{k\in \cI} |\psi_k'(w)|^\Theta_1\geq 2\quad \text{for every} \quad w\in U$$
 which itself implies that the topological 
pressure at $t=\Theta$ of $(\psi_k)_{k\in \cI}$ is strictly positif. Again by 
\cite{MauldinUrb02} or \cite{PUbook} follows finally that the Hausdorff dimension of the limit set 
of $(\psi_k)_{k\in \cI}$ is strictly larger than $\Theta$.
\epf


\medskip


\bibliographystyle{plain}



\end{document}